\documentclass[12pt,a4paper]{article}
\usepackage{amssymb,amsfonts,amsmath,amsthm,graphicx}
\makeindex

\theoremstyle{plain}
\newtheorem{teo}{Theorem}[section]
\newtheorem{defn}[teo]{Definition}
\newtheorem{lema}[teo]{Lemma}
\newtheorem{cor}[teo]{Corollary}
\newtheorem{prop}[teo]{Proposition}
\newtheorem{ex}[teo]{Example}
\newtheorem{obs}[teo]{Remark}

\newcommand{\sq}{\hfill\rule{2mm}{2mm}}

\newcommand{\2}{\vspace{.2cm}}

\begin{document}

\pagestyle{plain}
\begin{center}\large{\bf Partial Actions on Categories}\footnote{The second named author and the third named author were partially supported by Conselho Nacional de Desenvolvimento Cient\'{\i}fico
e Tecnol\'{o}gico (CNPq, Brazil) and  he second named author was also partially supported by FAPESP, (Funda\c c\~ao de Amparo a Pesquisa do Estado de S\~ao Paulo.).
}\end{center}

\begin{center}
{\bf {\rm Wagner Cortes$^1$, Miguel Ferrero$^2$, Eduardo N. Marcos$^3$}}\end{center}

\begin{center}{ \footnotesize 1,2 Instituto de Matem\' atica\\
Universidade Federal do Rio Grande do Sul\\
91509-900, Porto Alegre, RS, Brazil\\
3, Departamento de Matematica\\
IME-USP, Caixa Postal 66281,\\
05315-970, S\~ao Paulo-SP, Brazil\\
E-mail: {\it cortes@mat.ufrgs.br}, {\it mferrero@mat.ufrgs.br}, {\it enmarcos@ime.usp.br}}\end{center}

\begin{abstract} In this paper we introduce the definition of partial action on small $k$-categories
generalizing the similar well known notion of partial actions on
algebras. The point of view of partial action which we use in this
paper is the one which was introduced by Exel in his work on
$C^*$-algebras, see \cite{E}. Various generalizations were done
afterward, see \cite{CJ, DEP, DE, DFP}. Also we define the notion
partial skew category. We prove similar results to the ones in
\cite{CM}. Finally we show a result given conditions for a partial
action to have a globalization.
\end{abstract}

\section{Introduction}

\hspace{.5cm} The point of view of partial group actions which we
consider here was introduced in the context of operator algebras by
R. Exel in \cite{E}. Partial group actions are natural to be
consider from distinct points of view. A different way of looking at
it, of the one considered here, appeared earlier in \cite{GM}. From
the point of view considered in this paper a purely algebraic
treatment was given recently in \cite{DEP}, \cite{DE}. In
particular, several aspects of Galois theory can be generalized to
partial group actions, see \cite{DFP} (at least under the additional
assumption that the associated ideals are generated by central
idempotents). Recently, Caeneepel and Janssen in \cite{CJ} developed
the theory of partial (co)-action of a Hopf algebra and generalized
the Hopf Galois theory to this situation. Moreover,  other authors worked with partial actions of Hopf algebras, see  \cite{ME}.

In this paper we will extend the notion of partial actions to
categories. Recall that a category $C$ is said to be small if the
objects of $C$ is a set $C_0$, called the set of objects. Actually
we will need the following notion.

\begin{defn} Let $k$ be a commutative ring. A small weak not necessarily
associative $k$-category $D$ (WNNA, for short) consists of:

(1) A set $D_0$, called the set of objects of $D$;

(2) For each pair $(x,y)$ of objects of $D$ a $k$-module
$_{y}D_x$, called the set of morphisms from $x$ to $y$;

(3) For each $x,y,z$ in $D_0$, a $k$-bilinear map $\circ:
{_{z}D_y\times {_{y}D_x}}\to {_{z}D_x}$ called a composition,
where $\circ(f,g)$ will be denoted by $f\circ g$ ($\circ$ is not
necessarily associative). \end{defn}

\vspace{.2cm}

If the composition is associative, then $D$ is called a weak
category. If for every $x\in D_0$ there exists $1_x\in{ _{x}D_x}$
such that $f\circ 1_{x}=f$ and $1_{x}\circ g=g$, for every $f\in{
_{y}D_{x}}$ and $g\in { _{x}D_{y}}$, then we  say that $D$ is WNNA
category with identities or, more precisely, a not necessarily
associative category (i.e, a NNA category). A $k$-category is a WNNA
category with identities and associative composition.

The notions of subcategory, functors, full, faithful and so on,
can be defined for WNNA categories in a similar way that are
defined for categories.

Throughout this paper $k$ will always denote a commutative ring and
$C$  a small weak $k$-category. Actually the notion of weak
categories and not necessarily categories are not used strongly on
the paper. We consider weak categories because as in the algebras
situation the objects in our ideals do not need to have identities,
but once more the important case is when this happens. This is the
important case when we can show that there is an envelopping action.
So one very important case which we consider is the case where all
the ideals have local identities and in this case we can restrict
ourselves to $k$-categories.

Actions of groups on a small $k$-category $C$ were extensively
studied by several authors, see \cite{CM}, \cite{G} and the
references quoted  therein. The following question arises: is it
possible to generalize group actions on $k$-categories to the
partial situation? The  main purpose of this paper is to give a
positive answer to this question.

Let $G$ be a group. Recall that an action of $G$ on a small
$k$-category $C$ is  an action of $G$ on the set $C_0$ of objects of
$C$ and  a family of $k$-module isomorphisms
$s:{_{y}C_{x}}\rightarrow {_{sy}C_{sx}}$, for each $s\in G$ and for
each couple of objects $x$ and $y$ in $C_0$ and we have that
$s(gf)=(sg)(sf)$ in case $g$ and $f$ are morphisms which can be
composed in the category. Moreover, for elements $t,s\in G$ and a
morphism $f$ we have $(ts)f=t(sf)$ and $ef=f$, where $e$ is the
identity of $G$. A category $C$ together with an action of $G$ on
$C$ is called a $G$-category. We remark here that this notion can be
defined for WNNA categories in the same way.

Recalling the definition of an ideal in a category we give the
following

\begin{defn}
An ideal in a WNNA category $C$ is a collection $I$ of morphisms
such that if $f$ is in $I$ then $(gf)h$ and $g(fh)$ are in $I$
whenever $(gf)h$ and $g(fh)$ are defined. Moreover, if $C$ is a WNNA
$k$-category, for $I$ to be an ideal we require in addition  that
${_{a}I_{b}}$ is a $k$-submodule of  the $k$-module  ${_{a}C_{b}}$,
where ${_{a}I_{b}}$ denotes the set of all morphisms in
${_{a}C_{b}}$ which belong to $I$.
\end{defn}

Every ideal $I$ in a WNNA $k$-category $C$ can be looked as a WNNA
subcategory, also denoted by $I$. In this case $I_0 =C_0$ and
${_a}I_{b} = {_a}C_{b} \cap I,$ for any $a,b$ in $I_0$.

\vspace{.2cm}

In Section 2 we recall the definition of a partial action of a group
on a set $X$ and on a $k$-algebra $A$, and the partial skew group
algebra introduced by Exel and Dokuchaev \cite{DE}.

We also introduce the notion of partial orbit and show that the
family of partial orbits of a set $X$ form a partition of it,
which is a generalization of what happens in the case of global
actions.

In Section 3 we introduce the notion of a partial action of a group
on a weak $k$-categories and  define the partial skew category. We
prove a coherence result between our approach and the
ring-theoretical approach, in the case the weak category has a
finite number of objects. Moreover, we show that the partial skew
category is equivalent to the full subcategory of the partial skew
category formed by taking one element in each equivalent class.

In Section 4 we give conditions for a partial action of a group on a
small $k$-category $C$ to have an enveloping action.

\begin{defn}
 Let $C$ be a $k$-WNNA  category, $x$ an object of $C$ and $I$ an ideal of $C$. A morphism $e$  in ${_{x}I_x}$ is called a local identity if, $e$ is an idempotent, $ef=f$
for all $f\in {_{x}I_{y}}$, and $fe=f$ for all $f\in {_{y}I_{x}}$. Morever the local identity is called central if $fe=ef$ for all $f\in {_{x}C_{x}}$.\end{defn}

 It is convenient to point out that if $C$ is a $k$-category and $I$ is an ideal of $C$ such that  for each object $x\in C$ there is  a morphism $e\in {_{x}I_{x}}$ that is a local identity, then $I$ is itself a small $k$-category.

\section{Partial skew group algebras}

Let $G$ be a group and $X$ a set. A {\it partial
 action} $\alpha$ of $G$ on $X$ is a collection of subsets
 $S_g$, $g\in G$, of $X$ and bijections
  $\alpha_g:S_{g^{-1}}\to S_g$ such
 that:

\2

 (i) $S_1=X$ and $\alpha_1$ is the identity mapping of $X$;

\2

 (ii) $S_{(gh)^{-1}}\supseteq \alpha_h^{-1}(S_h\cap S_{g^{-1}})$;

\2

(iii) $\alpha_{g}\circ \alpha_{h}(x)=\alpha_{gh}(x)$, for any $x\in
\alpha^{-1}_{h}(S_{h}\cap S_{g^{-1}})$.

\2

{\bf Remark:}

1.  The property (ii) is equivalent to $\alpha_{g}(S_{g^{-1}}\cap S_{h})=S_{g}\cap S_{gh}$,
 for all $g,h\in G$.

2. We also have that  $\alpha_{g^{-1}}=\alpha^{-1}_{g}$, for every $g\in
 G$.

Following \cite{A},  given a partial action $\alpha$ of $G$ on $X$,
a globalization of $\alpha$, also called an enveloping action,  is a
pair $(Y,\beta)$ such that $X$ can be considered as a subset of $Y$,
$\beta$ is a global action of $G$ on $Y$, $Y=\cup_{g\in
G}\beta_g(X)$, $S_g=X\cap \beta_g(X)$ and $\alpha_g:S_{g^{-1}}\to
S_g$ is equal to $\beta_g|_{S_{g^{-1}}}$, $g\in G$. In other words,
$\alpha$ is the restriction of $\beta$ to $X$.

If we have an additional structure some conditions can be imposed to
the subsets $S_g$ and the maps $\alpha_g$. For example, for partial
actions on topological spaces all the $S_g$ are open subspaces of
$X$ and all the mappings $\alpha_g$ are homeomorphisms of
topological spaces. In Theorem 1.1 of \cite{A} the author proved
that globalization does exist for partial actions on topological
spaces, which clearly implies the result for partial actions on
sets.

For the definition of a partial action of a group $G$ on a
$k$-algebra the authors considered that any $S_g$ is an ideal
of $R$ and that every $\alpha_g$ is an isomorphism of algebras, $g\in G$.

Let $\alpha$ be a partial action of $G$ on the $k$-algebra $R$.
The partial skew group algebra $R\star_{\alpha}G$, see \cite{DE}, is
defined as the set of all finite formal sums $\sum_{g\in
G}a_{g}u_{g}$, $a_g\in S_g$ for every $g\in G$, where the addition
is defined in the usual way and the multiplication is determined
by
$(a_{g}u_{g})(b_{h}u_{h})=\alpha_{g}(\alpha_{g^{-1}}(a_{g})b_{h})u_{gh}$.

For the sake of completeness we recall now some facts and give a
proposition.

Let $A$ be a $k$-algebra with identity element and $\{e_i\}_{i=1}^n$
a complete set of orthogonal idempotents, i.e., a set of orthogonal
idempotents whose sum is the identity of $A$. Then we define a
$k$-category with a finite number of objects, denoted either by
$(C_A, \{e_i\}_{i=0}^{n})$ or $C_A$, as follows:

$Obj(C_A) =\{e_i: 1\leq i\leq n\}$ and $Hom(e_i, e_j) = e_i A
e_j$, for all $1\leq i, j\leq n$. Composition is defined in a natural way via
the product $(e_iAe_j)(e_jAe_k) \subset e_iAe_k$.

Conversely, given a $k$-category $C$ which has a finite number of
objects $\{e_i: 1\leq i\leq n\}$, we define the, so called,
$k$-algebra of homomorphism $a(C)$ in the following way: $a(C)$,
as a set, is equal to $\oplus_{i,j}Hom(e_i,e_j)$. Note that the
elements of $a(C)$ can be seen as  matrices. So the addition and
multiplication in $a(C)$ is defined as for matrices.

The following proposition is probably known nevertheless we give a proof here.

\begin{prop}\label{21}
(i) Let $A$ be a $k$-algebra with identity element and
$\{e_i\}_{i=1}^n$ a complete set of  orthogonal idempotents of $A$.
Then the associated $k$-algebra $a(C_{A})$ of $C_{A}$ is isomorphic
to $A$.

(ii) Let $C$ be a $k$-category with a finite number of objects.
Then the categories $C$ and $C_{a(C)}$ are equivalent.
\end{prop}

{\bf Proof.} (i) We define $\varphi:A\rightarrow a(C_{A})$ by
$\varphi(r)=(e_{i}re_{j})_{i,j}$ and $\psi:a(C_{A}) \rightarrow A$
by $\psi((e_{i}r_{i,j} e_{j})_{i,j})=\sum_{i,j}e_{i}r_{i,j}e_{j}$.
It is easy to see that $\varphi$ is a homomorphism of algebras,
$\varphi\circ\psi=id_{a(C_A)}$ and $\psi\circ\varphi=id_A$.

(ii) Denote by $\{a_1,...,a_n\}$ the objects of the category $C$. By
definition $a(C)=\oplus{_{i,j}}Hom(a_i,a_j)$ with matrix operations.
We define the category $C_{a(C)}$ by taking the complete set of
orthogonal idempotents $\{_{a_1}1_{a_1},...,{_{a_n}1_{a_n}}\}$ in
$a(C)$. Let $T$ be the functor from $C$ into $C_{a(C)}$ defined by
$T(a_i)={_{a_i}}1_{a_i}$, for any $1\leq i\leq n$, and for each
$f:a_i\to a_j$, $T(f)={_{a_i}}1_{a_i}f{_{a_j}}1_{a_j}$. It can
easily be seen that $T$ is an equivalence of categories. \sq

\begin{defn} Let $\alpha$ be a partial action of $G$ on a set $X$.  For
each $a,b\in X$ we say that $a,b$ are $\alpha$-equivalent,
$a{\sim_{\alpha}} b$, if there exists $g\in G$ such that $a\in
S_{g^{-1}}$ and $b=\alpha_{g}(a)$. Briefly $a{\sim_{\alpha}} b$
will be denoted by $a\sim b$.
\end{defn}

\begin{lema} Let $\alpha$ be a partial action of $G$ on a set $X$.
Then the relation $\sim$ is an equivalence
relation.\end{lema}

Proof. Straightforward. \sq

\vspace{.4cm}

Let $\alpha$ be a partial action of a group $G$ on a set $X$. For
each $x\in X$ the partial orbit of $x$ is
$H^{\alpha}(x)=\{\alpha_{g}(x):x\in S_{g^{-1}}\}$. By Lemma 2.3
the set of all partial orbits form a partition of $X$.

\begin{ex} \label{e1} {\rm Consider $X=\{e_{1},e_{2},e_{3},e_{4}\}$
and denote by $G$ the cyclic group generated by $\sigma$ of order
5. Let us take the subsets $S_{1}=X, S_{\sigma}=\{e_{2}\},
S_{\sigma^{2}}=\{e_{4}\},S_{\sigma^{3}}=\{e_{3}\},S_{\sigma^{4}}=\{e_{1}\}$
and define $\alpha$ by $\alpha_{1}=id_{X}$,
$\alpha_{\sigma}(e_{1})=e_{2}$,
$\alpha_{\sigma^{2}}(e_{3})=e_{4}$,
$\alpha_{\sigma^{3}}(e_{4})=e_{3}$,
$\alpha_{\sigma^{4}}(e_{2})=e_{1}$. It is easy to see that
$\alpha$ is a partial action of $G$ on $X$,
$H^{\alpha}(e_{1})=\{e_{1},e_{2}\}=H^{\alpha}(e_{2})$ and
$H^{\alpha}(e_{3})=\{e_{3},e_{4}\}=H^{\alpha}(e_{4})$}.\end{ex}

Assume that $\alpha$ is a partial action of $G$ on $X$ and let
$(Y, \beta)$ be the globalization of $\alpha$. Thus for any $g\in
G$, $\beta_g$ is a bijection of $Y$. Hence $\beta$ is an action on
$Y$ and so it defines an equivalence relation denoted by $\approx
$: for $u,v\in Y$, $u\approx v$ if there exists $g\in G$ such that
$\beta_g(u)=v$. We denote the orbit of $y\in Y$ by $H^{\beta}(y)$.

We have the following

\begin{lema} Let $\alpha$ be a partial action of $G$ on $X$ and
let $(Y,\beta)$ be its enveloping action. Then the equivalence
relation defined by $\alpha$ on $X$ is the restriction of the
equivalence relation defined by $\beta$ on $Y$. In particular, for
any $x\in X$ $H^{\alpha}(x)=H^{\beta}(x)\cap X$. \end{lema}

Proof. Straightforward. \sq

\vspace{.4cm}

The following remark will be useful.

\begin{obs} $($See \cite{CM}$)$ Let $C$ be a $k$-category with a finite number of
objects $\{f_{i}: 1\leq i\leq n\}$. For each $i\in \{1,...,n\}$ we
associate the projective $C$-module $C_{f{i}}(f_j) =
{_{f_{j}}C_{f_{i}}}$. Note that $\sum_{1\leq i\leq n} \oplus
C_{f_i}\simeq a(C)$. Hence
$End_C(\coprod_{i=1}^{n}C_{f_{i}})=End(a(C))\simeq a(C)$.
\end{obs}

\begin{teo} Let $G$ be a finite group and $\alpha$ a partial action of
$G$ on a $k$-algebra with identity $A$.  Suppose that
$S_{g}=Ae_{g}$, where $(e_g)_{g\in G}$ is a set of orthogonal
idempotents whose sum is 1 such that  $Ae_g=e_gA$, for all $g\in G$.
Assume that $\{H_{1},...,H_{n}\}$ is the family of partial orbits of
the set $\{e_{g}:g\in G\}$. Let $\{f_1,\cdots f_n\}$ be the
idempotents, chosen one for each orbit. Then
$End_{A*_{\alpha}G}(\coprod_{i=1}^{n}A*_{\alpha}Gf_{i})$ is Morita
equivalent to $A*_{\alpha}G$. In particular, if the idempotents are
all in the orbit of a fix idempotent $e$
 then  $End_{A*_{\alpha}G}(A*_{\alpha}Ge)$ is Morita
equivalent to $A*_{\alpha}G$.\end{teo}

{\bf Proof.} Note that the ideals $S_{g}$ are idempotents and by
(\cite{DE}, Theorem 3.1), the partial skew group algebra
$A*_{\alpha}G$ is an associative algebra. By Proposition \ref{21}
$A*_{\alpha}G$ is isomorphic to $a(C_{A*_{\alpha}G})$ and
$C_{A*_{\alpha}G}$ is isomorphic to the full subcategory formed by
the objects $\{f_{1},...,f_{n}\}$, where $f_{i}\in H_{i}$. Hence,
$a(C_{A*_{\alpha}G})$ is Morita equivalent to
$a(C\{f_{i}\}_{i=1}^{n})$. So, by Remark 2.6
$\coprod_{i,j=1}^{n}f_{j}A*_{\alpha}Gf_{i}$ is isomorphic to
$End_{A*_{\alpha}G}(\coprod_{i=1}^{n}A*_{\alpha}Gf_{i})$.\sq

\begin{ex} As in Example \ref{e1}, take
$S=\bigoplus_{i=1}^{4}Re_{i}$, where $R$ is a commutative ring, with
the same partial action as before. Using the above theorem we get
that $T=S*_{\alpha}G$ is Morita equivalent to
$End(e_{1}T+e_{2}T)$\end{ex}

\section{Partial actions of groups on weak categories}

\hspace{.5cm} In the remaining of the paper $C$ will denote a small
weak $k$-category, where $k$ is commutative ring and $C_0$ the set
of objects of $C$.

In the next definition we will give first a partial action
$\alpha_0$ of  a group $G$ on the set $C_0$. So for any $g\in G$ a
subset $C_0^g$ is given and $\alpha_0^g:C_0^{g^{-1}}\to C_0^g$ is a
bijection. If $x\in C_0^{g^{-1}}$, $\alpha_0^g(x)$ will be denoted
by $gx$.

\begin{defn} Let $G$ be a group. We say that $\alpha=\{\alpha^g | g\in G\}$ is a partial action of $G$
on $C$ if the following conditions hold:

(i) $G$ acts partially on the set of objects $C_0$ of $C$. This
partial action will be denoted by $\alpha_0$ and the subsets
associated to this partial action by $C_{0}^{g}$, $g\in G$;

(ii) For each $g\in G$ there exists an ideal  $\mathcal{I}^{g}$ of
$C$ such that $_{a}\mathcal{I}_{b}^{g} = 0$ if one of the elements
$a$ or $b$ are not in $C_0^g$;

(iii) There are equivalence of weak categories
$\alpha^{g}:\mathcal{I}^{g^{-1}}\rightarrow\mathcal{I}^{g}$, for any
$g\in G$, such that for $f\in {_y\mathcal{I}}_x^{g^{-1}}$,
$\alpha^g(f)\in {_{gy}\mathcal{I}}_{gx}^{g}$, where $x,y$ are in
$C_0^{g^{-1}}$;

(iv) $\mathcal{I}^e = C$ and $\alpha^e = Id$;

(v) For any pair of objects $(x,y) \in C_0\times C_0$ we have that
\begin{center} $\alpha^{h^{-1}}(_{y}\mathcal{I}_{x}^{h} \cap
{_{y}\mathcal{I}}_{x}^{g^{-1}})\subseteq
{_{h^{-1}y}\mathcal{I}}_{h^{-1}x}^{(gh)^{-1}}$, if $x,y\in
C_{0}^{h}$;\end{center}

(vi) If  $x,y\in C_0^{h}\cap C_{0}^{g^{-1}}$ and $f\in
\alpha^{h^{-1}}(_{y}\mathcal{I}_{x}^{h} \cap
{_{y}\mathcal{I}}_{x}^{g^{-1}})$, then
$\alpha^{g}(\alpha^{h}(f))=\alpha^{gh}(f)$.\end{defn}

Note that the conditions (v) and (vi) above fit with the conditions
(ii) and (iii) of the definition of partial actions of groups on
algebras given in \cite{DE}. Also, as in \cite{DE}, it can easily be
checked that condition (iv) can be replaced by the condition
$\alpha^{h^{-1}}(_{y}\mathcal{I}_{x}^{h} \cap
{_{y}\mathcal{I}}_{x}^{g^{-1}})={_{h^{-1}y}\mathcal{I}}_{h^{-1}x}^{(gh)^{-1}}\cap
{_{h^{-1}y}\mathcal{I}}_{h^{-1}x}^{h^{-1}}$.

\vspace{.2cm} Now we give natural examples of partial actions on
small weak $k$-categories.

\begin{ex} {\rm Assume that $C$ is a  small weak $k$-category and $\beta$ is a global
action of $G$ on $C$. Let $D$ be an ideal of $C$. We define a
partial action $\alpha$ of $G$ on $D$ by restriction of $\beta$ to
$D$ as follows:

The partial action is a global action on $D_0=C_0$

The ideals $\mathcal{I}^g$ of $D$ are defined by
$_y\mathcal{I}_x^g={_yD_x}\cap \beta_g(_{g^{-1}y}D_{g^{-1}x})$,
where, as above, $\beta_{g^{-1}}(x)$ is denoted simply by $g^{-1}x$.

Let us show first that $\mathcal{I}^g$ is an ideal of $D$. Let
$f\in {_y\mathcal{I}}_x^g$ and $l, m$ morphisms such that $l\in
{_zD_y}$ and $m\in{_xD_u}$, with $u$ and $z$ in $C_0=D_{0}$. Then
there are \begin{center}$t\in {_{g^{-1}y}D}_{g^{-1}x}$,
$\tilde{l}\in {_{g^{-1}z}C}_{g^{-1}y}$, $\tilde{m}\in
{_{g^{-1}x}C}_{g^{-1}u}$
\end{center}  with
$f=\beta_g(t)$, $l=\beta_g(\tilde{l})$ and $m=\beta_g(\tilde{m})$.
So
$lfm=\beta_g(\tilde{l})\beta_g(t)\beta_g(\tilde{m})=\beta_g(\tilde{l}t\tilde{m})\in
{_zD_u}\cap \beta_g({_{g^{-1}z}D}_{g^{-1}u})
={_z\mathcal{I}}_u^g$.

It is easy to check that $\alpha^g:{_y\mathcal{I}_x^{g^{-1}}}\to
{_{gy}\mathcal{I}}^g_{gx}$ is a bijection and an isomorphism of
$k$-modules, for any $x,y\in D_0^g$.

Now let $f\in \alpha^{h^{-1}}(_y\mathcal{I}_x^{h} \cap
{_y\mathcal{I}_x^{g^{-1}}})$, i.e., $$\alpha^h(f)\in
{_{y}\mathcal{I}_{x}^{h}} \cap
{_{y}\mathcal{I}_{x}^{g^{-1}}}={_yD_x\cap
\beta_h(_{h^{-1}y}D_{h^{-1}x}})\cap \beta_{g^{-1}}(_{gy}D_{gx}).$$
Hence $f\in
\beta_{h^{-1}}(\beta_{g^{-1}}(_{gy}D_{gx}))=\beta_{(gh)^{-1}}(_{gy}D_{gx})$
and consequently $$f\in
 {_{h^{-1}y}D_{h^{-1}x}}\cap \beta_{(gh)^{-1}}(_{gy}D_{gx})=
 {_{h^{-1}y}\mathcal{I}}_{h^{-1}x}^{(gh)^{-1}}.$$ Thus condition
 $(iv)$
 of the definition of partial actions is satisfied. Finally,
 condition (v) also holds since $\alpha^g$ is defined as
 restriction of $\beta_g$}. \sq
\end{ex}

\begin{obs} {\rm Note that we can give a more general example if we
change slightly the definition of ideal of a weak category. In the
following an ideal $D$ of $C$ will be a subset $D_0$ of $C_0$
together with a collection of morphisms ${_y}D_x$, for every $x,y\in
D_0$, satisfying the same conditions of the Definition 1.1}.

\end{obs}

\begin{ex} {\rm Assume that $C$ is a  small weak $k$-category and
$\beta$ is a global action of $G$ on $C$. Let $D$ be an ideal of $C$
in the sense of the above remark. For any $g\in G$ we put
$D_0^g=D_0\cap \beta_g(D_0)$. The partial action $\alpha_0$ on $D_0$
is defined as the restriction of $\beta$, i.e.
$\alpha_0^g:D_0^{g^{-1}}\to D_0^g$ is equal to
$\beta|_{D_0^{g^{-1}}}$. The rest is defined as in the above
example. The ideals $\mathcal{I}^g$ of $D$ are defined by
$_y\mathcal{I}_x^g={_yD_x}\cap \beta_g(_{g^{-1}y}D_{g^{-1}x})$ and
$\alpha^g$ is the restriction of $\beta_g$ to $\mathcal{I}^{g^{-1}}$
It is easy to show, as above, that this gives a partial action of
$G$ on $D$ which not global in $D_0$.}

\end{ex}

Now we introduce a partial version of skew category.

\begin{defn} Let $\alpha$ be partial action of a group $G$ on a small
weak $k$-category $C$. We define the skew WNNA category
$C*_{\alpha}G$ as follows:

(i) $(C*_{\alpha}G)_{0}=C_{0}$.

(ii) For each $x,y\in C_{0}$, $_{y}(C*_{\alpha}G)_{x}=\oplus_{g\in
G}\ _{y}\mathcal{I}_{gx}^{g}$.

For each $f\in \ _{z}\mathcal{I}_{ty}^{t}$, $l\in\
_{y}\mathcal{I}_{gx}^{g}$  we define the composition  by the
following rule: $fl=\alpha^{t}(\alpha^{t^{-1}}(f)\circ l)\in
{_{z}\mathcal{I}_{(tg)x}^{tg}}$.\end{defn}

\begin{ex} Let $G$ be a group and $\alpha$ a partial action
on a small $k$-category $C$. If for any $g\in G$ the ideals
$\mathcal{I}^{g}=C$, then for all $x,y\in C_{0}$,
$\alpha_g:{_{x}C_{y}}\rightarrow {_{gx}C_{gy}}$ is an isomorphism
of $k$-modules. Note that the action of $G$ on $C_0$ is global in this case. Thus, $\alpha$ is a global action of $G$ on $C$ and
$C*_{\alpha}G$ defined above is the ordinary skew category $C[G]$,
see \cite{CM}.
\end{ex}

We know that the partial skew group algebra introduced by Dokuchaev
and Exel is not necessarily associative (Example 3.5 of \cite{DE}).
Similarly the composition map in $C*_{\alpha}G$ is not necessarily
associative, in general. The case where it is associative is of
special interest, because of this we give the next definition.

\begin{defn} Let $G$ be a group and $\alpha$ a partial action of G on
a small $k$-category $C$. We say that the partial action $\alpha$ is
associative if the composition of maps in $C*_{\alpha}G$ is
associative.
\end{defn}

\begin{obs} As a consequence of the definition above, if $C$ is small $k$-category
and the partial action $\alpha$ is associative, then $C*_{\alpha}G$
is a category and we call it the partial skew category.\end{obs}

In the next theorem we assume that $\alpha$ is a partial action
of $G$ on a small $k$-category $C$ with finite number of objects  such that $\alpha$ is not associative. In this case, $a(C*_{\alpha}G)$ is not necessary associative $k$-algebra.

\begin{teo} Let $\alpha$ be a partial action of a group $G$ on a
small $k$-category $C$ with a finite number of objects. Then $G$ acts partially on $a(C)=
\oplus_{x,y\in C_{0}}\ _{y}C_{x}$ and $a(C*_{\alpha}G)$ is
isomorphic to $a(C)*_{\alpha} G$.\end{teo}

 Proof. For each $g\in
G$ let $a(C)_g=\oplus _{x,y\in C_{0}}\ _{y}\mathcal{I}_{x}^{g}$ is
an ideal of $a(C)$ and $\alpha_g:a(C)_{g^{-1}}\rightarrow a(C)_g$,
defined by
$\alpha_g|_{{_{y}\mathcal{I}^{g^{-1}}_{x}}}=\alpha^g|_{{_{y}\mathcal{I}^{g^{-1}}_{x}}}$,
for all $x,y\in C_0^{g^{-1}}$, is an isomorphism of ideals.

Now we show that $\alpha$ is a partial action of $G$ on $a(C)$.
The first condition of the definition of partial actions is
obvious. For the second suppose that $f\in
\alpha_{h^{-1}}(a(C)_h\cap a(C)_{g^{-1}})$. We can assume that
$f\in {_{y}C_x}$ , so $\alpha_{h}(f)\in {_{hy}C_{hx}}$  and consequently $f\in a(C)_{(gh)^{-1}}$. Thus
the second condition of the definition of partial actions is
fulfilled. Finally the condition (iii) of Definition 1.1 in
\cite{DE} follows immediately from condition (vi) of Definition
3.1.

We define $\varphi: a(C*_{\alpha}G)\rightarrow a(C)*_{\alpha}G$ by
$\varphi(f_{g})=f_gu_g$, where $f_g$ is an elementary morphism in
$_{y}\mathcal{I}^g_{gx}\subseteq{_{y}(C*_{\alpha} G)_x}$ and
$\{u_g\mid g\in G\}$ denotes the  canonical generators of $a(C)\star_{\alpha}
G$. We clearly have that $\varphi$ is a well defined homomorphism
of $k$-algebras. Finally $\Psi: a(C)*_{\alpha}G\rightarrow
a(C*_{\alpha}G)$ defined by $\Psi(f_g u_g)=f_{g}$, for any $f_g\in
a(C)_g$, is clearly an inverse of $\varphi$. \sq

\vspace{.4cm}

Recall that an algebra $A$ is strongly associative if for any
partial action $\alpha$ of a group $G$ on $A$ is always
associative. A semiprime algebra is strongly associative
(\cite{DE}, Corollary 3.4)

The following is immediate from Theorem 3.7.

\begin{cor} Let $G$ be a group, $A$ a strongly associative $k$-algebra and
$C_{A}$ the category with a single object $\{A\}$ and endomorphism
$k$-algebra $A$. Suppose that $G$ acts partially on $C_{A}$. Then
$a(C_{A}*_{\alpha}G)\simeq A*_{\alpha}G$ and so $\alpha$ is
associative. In particular, if $A$ is semiprime then the
associated skew NNA category is associative for any partial action
$\alpha$ of a group $G$ on $A$ .\end{cor}

Let $\alpha$ be a partial action of $G$ on a small $k$-category $C$.
We define the following relation: $x\sim y$ if there exists $g\in G$
such that $_{x}1_x\in {_{x}\mathcal{I}_{x}^{g^{-1}}}$ and
$y=\alpha^g(x)$, where $1_x\in {_{x}\mathcal{I}_{x}^{g^{-1}}}$
denotes the identity morphism from $x$ to $x$.

The proof of following lemma is standard.

\begin{lema} The relation $\sim$ is an equivalence relation.\end{lema}

The next proposition has a similar proof as Lemma 2.5, of
\cite{CM}. For the sake of completeness we give a proof here,
adapted to our case.

\begin{prop} Let $\alpha$ be a partial action of a group $G$ on a small $k$-category
$C$ such that $\alpha$ is associative. If $x\sim y$, then the
objects $x$ and $y$ are isomorphic in $C*_{\alpha}G$.
\end{prop}

Proof. Suppose that $x$, $y$ are equivalent, i.e, there exists $h\in
G$ such that $y=hx$ and $_{x}1_x\in {_{x}\mathcal{I}_{x}^{g^{-1}}}$.
Since $_{y}1_y\in {_{y}\mathcal{I}^{h}_{y}}$, we have that
$_{y}1_{y}={_{y}1_{hx}}\in {_{y}\mathcal{I}^{h}_{hx}}\subseteq
{_{y}(C*_{\alpha}G)_{x}}$. On the other hand,
$_{x}1_{x}={_{x}1_{h^{-1}y}}\in
{_{x}\mathcal{I}^{h^{-1}}_{h^{-1}y}}\subset {_{x}(C*_{\alpha}G)_{y}}$.
We claim that $_{x}1_{x}={_{x}1_{h^{-1}y}}\circ _{y}1_{hx}$. In
fact, $_{x}1_{h^{-1}y}\circ
{_{y}1_{hx}}=\alpha^{h^{-1}}(\alpha^{h}(_{x}1_{h^{-1}y})\circ
{_{y}1_{hx}})=\alpha^{h^{-1}}(_{hx}1_{y}\circ
{_{y}1_{hx}})=\alpha^{h^{-1}}(_{hx}1_{hx})=\!\!_{x}1_{x}$. Using
similar methods we can show that
$_{y}1_{y}={_{y}1_{hx}}\circ{_{x}1_{h^{-1}y}}$.\sq

\vspace{.4cm}

The next corollary is a direct consequence of the above result.

\begin{cor} Let $\alpha$ be a partial action of a group $G$ on a small
$k$-category $C$ such that $\alpha$ is  associative and  $S$ a representative set of the equivalence relation defined before (that is, there is in $S$ exactly one element of each equivalence class)
and $C^*$ be the full subcategory of $C*_{\alpha}G$ whose objects are the elements of $S$. Then $C^*$ is
equivalent to $C*_{\alpha}G$ .\end{cor}

\section{Globalization of partial actions}

\hspace{.5cm} In this section we consider always small
$k$-categories, unless otherwise stated.

Examples 3.2 and 3.4 are natural examples of partial actions of
groups on small $k$-categories which can be obtained by restriction
of global actions to ideals. Moreover, if a partial action $\alpha$
is obtained in that way, then $C*_{\alpha}G$, defined as before, is
an associative category. Thus it is natural to ask when a partial
action can be obtained by restriction of a global action. This
question has been considered in \cite{A} for partial actions on
topological spaces and in \cite{DE} for partial action on algebras
with identity element.

\begin{defn} Let $B$ and $D$ be two WNNA categories. We say that $T:B\rightarrow D$
is a quasi-functor if the following conditions are satisfied:

(i) For each object $b$ of $B$, $T(b)$ is an object of $D$;

(ii) For each morphism $f$ in $B$, $T(f)$ is a morphism in $D$;

(iii) Given two morphims $f$,$g$ in $B$ such that $\exists f\circ g$
in $B$ we have that $T(f)\circ T(g)$ does exist in $D$ and $T(f\circ
g)=T(f)\circ T(g)$.\end{defn}

Note that if $B$ and $D$ are categories, then a functor $T:B\to D$
is a quasi-functor such that $T(_{x}1_x)={_{T(x)}}1_{T(x)}$, for all
$x\in B_0$.

Induced by the definition given in \cite{DE}, Section 4, we give
the following.

\begin{defn}\label{defi} Let $(C,\alpha)$ be a small $k$-category $C$ together
with a partial action $\alpha$ of $G$ on $C$. We say that a pair
$(D,\beta)$, where $D$ is a $k$-category and $\beta$ is a global
action of $G$ on $D$, is an enveloping (also called a globalization)
of $(C,\alpha)$ if the following conditions are satisfied:

(i) There is a faithful quasi-functor $j:C\to D$;

(ii) For each $f\in {_{y}j(C)_{x}}$, $g\in {_{z}D_{y}}$ and $h\in
{_{x}D_{v}}$, where $x,y,z,v\in j(C_{0})$, we have $gfh\in
{_{z}j(C)_{v}}$;

(iii) $j(_y\mathcal{I}_x^g)=j(_yC_x)\cap
\beta_g(j(_{{g^{-1}}y}C_{g^{-1}x}))$, for all $x,y\in C_0$;

(iv) $j\circ \alpha^g(f)= \beta_g\circ j(f)$, for any $f\in
{_y\mathcal{I}^{g^{-1}}_x}$;

(v) $_{y}D_{x}=\sum_{g\in G} \beta_g(j(_{g^{-1}y}C_{g^{-1}x}))$, for
any $x,y\in D_0$ such that $g^{-1}x,g^{-1}y\in C_{0}$.

It is convenient to remark that when $(C,\alpha)$ has an enveloping action $(D,\beta)$ we have that $C*_{\alpha}G$ is a subcategory of the skew category $D[G]$.

\end{defn}

\begin{defn}\label{defi2} Given small $k$-categories $D$ and
$D^{'}$, we say that global actions $(D,\beta)$ and
$(D^{'},\beta^{'})$ of $G$ on $D$ and $D^{'}$, respectively, are
equivalent if there exists an equivalence of categories $\Phi:D\to
D^{'}$ such that for any $g\in G$ we have $\beta_g\circ
\Phi=\Phi\circ \beta_g^{'}$.
\end{defn}

\begin{lema} Assume that $\alpha$ is a partial action of a group
$G$ on a small $k$-category $C$ which has a globalization
$(D,\beta)$ and let $j:C\to D$ the canonical faithful functor.
Then for any $x\in C_0$ and $g_1,...,g_n\in G$ the submodule
$\sum_{1\leq i\leq n}\beta_{g_i}(j(_{{g_i^{-1}}x}C_{g_i^{-1}x}))$
of $_x{D_x}$ has an identity element with respect to composition.
\end{lema}

{\bf Proof.} By induction it is enough to prove the result for
$n=2$. Put
$N=\beta_g(j(_{g^{-1}x}C_{g^{-1}x}))+\beta_h(j(_{h^{-1}x}C_{h^{-1}x}))$.
Since $\beta_g(j(_{g^{-1}x}1_{g^{-1}x}))$ is an identity for
$\beta_g(j(_{g^{-1}x}C_{g^{-1}x}))$ and
$\beta_h(j(_{h^{-1}x}1_{h^{-1}x}))$ an identity for
$\beta_h(j(_{h^{-1}x}C_{h^{-1}x}))$, it is easy to see that
$$\beta_g(j(_{g^{-1}x}1_{g^{-1}x}))+\beta_h(j(_{h^{-1}x}1_{h^{-1}x}))-
\beta_g(j(_{g^{-1}x}1_{g^{-1}x}))\beta_h(j(_{h^{-1}x}1_{h^{-1}x}))$$
is an identity for $N$. \sq

\vspace{.4cm}

In the next result we will assume that the $k$-subspace
$_x\mathcal{I}_x^g$ contains a local identity, for any $g\in G$ and
$x\in C_0$. We should point out that this local identity is not
necessarily the identity in $_xC_x$. Now we prove the main
theorem of this section.

\begin{teo} \label{imp} Let $\alpha$ be a  partial action of a group
$G$ on a small $k$-category $C$ such that $\alpha_0$ is global on
$C_0$. Then there exists an enveloping action of $(C,\alpha)$ if and
only if all the $k$-spaces $_x\mathcal{I}_x^g$ contains a local identity
element, for any $x\in C_0$ and $g\in G$. Moreover, the enveloping
action, if does exists, it is unique up to equivalence.
\end{teo}

{\bf Proof.} Assume that $(D,\beta)$ is an enveloping action of
$(C,\alpha)$ and denote by $j:C\to D$ the  functor of Definition
4.2. Note that, in this case, since $\alpha_0$ is global in $C_0$ we
can assume that $C_0=D_0$ and so $j(C)$ is an ideal of $D$. Now
$_x\mathcal{I}_x^g$ has an identity
$j^{-1}(j(_x1_x)\beta_g(j(_{g^{-1}x}1_{g^{-1}x})))$, for any $g\in
G$ and $x\in C_0$, where $_x1_x$ denotes the identity of $_xC_x$.

Conversely, in the rest of the proof we assume that
$_x\mathcal{I}_x^g$ contains an identity $_x1_x^g$, for any $x\in
C_0$ and $g\in G$, that we shortly denote by $1^g$ when there is
no possibility of misunderstanding. Define the category $B$ as
follows: $B_0=C_0$ and for any $x,y\in B_0$ the $k$-module of
morphisms $_yB_x$ is defined as the $k$-module $F(G,{_yC_x})$ of
all the maps $\sigma$ from $G$ to the direct product $\prod_{g\in
G} {_{gy}C_{gx}}$ such that $\sigma(g)\in
{_{g^{-1}y}C_{g^{-1}x}}$, for any $g\in G$.

As in \cite{DE} we write $\sigma\mid_g$ to denote $\sigma(g)$. The
composition of the morphisms $\sigma\in {_yB_x}$ and $\tau\in
{_zB_y}$ is defined by $\tau\circ
\sigma\mid_h=\tau\mid_h\circ\sigma\mid_h$.

We define a global action $\beta$ of $G$ on $B$ as follows: if
$\sigma\in {_yB_x}$ and $h\in G$, then we put
$\beta_h(\sigma)\mid_g =\sigma\mid_{h^{-1}g}$. Since
$\sigma\mid_{h^{-1}g}\, \in \,_{g^{-1}hy}C_{g^{-1}hx}$ it follows
that $\beta_h(\sigma)\in {_{hy}B_{hx}}$ and so $\beta$ is
well-defined. It is easy to see that $\beta$ is an action of $G$
on $B$.

Now we define a functor $j:C\to B$. As map from $C_0$ to $B_0$,
$j$ is the identity. If $\sigma\in {_yC_x}$ we put
$j(\sigma)\mid_h=\alpha^{h^{-1}}(\sigma {1^h})\in
{_{h^{-1}y}C_{h^{-1}x}}$. Thus $j(\sigma)\in {_yB_x}$ and
consequently $j$ is well-defined. We see that $j$ is a faithful
functor from $C$ into $B$. Suppose that $\mu, \eta\in C$ and
$j(\mu)=j(\eta)$. Then for any $h\in G$ we have
$j(\mu)\mid_h=j(\eta)\mid_h$ and so $\alpha^{h^{-1}}(\mu
{1^h})=\alpha^{h^{-1}}(\eta {1^h})$. Taking $h=e$, the identity of
$G$, we obtain $\mu=\eta$.

Now assume that $\mu\in {_yC_x}$ and $\eta\in {_zC_y}$. Hence
$\eta \mu\in {_zC_x}$ and
$j(\eta\mu)\mid_h=\alpha^{h^{-1}}(\eta\mu1^h)$. On the other hand
\begin{center}$j(\eta)j(\mu)\mid_h=j(\eta)\mid_h\circ
j(\mu)\mid_h=\alpha^{h^{-1}}(\eta 1^h)\circ \alpha^{h^{-1}}(\mu
1^h)=\alpha^{h^{-1}}(\eta 1^h\mu 1^h)=\alpha^{h^{-1}}(\eta\mu
1^h)$.\end{center} Therefore $j(\eta\mu)=j(\eta)j(\mu)$.

Let $D$ be the subcategory of $B$ defined as follows: the set of
objects $D_0$ of $D$ is equal to $C_0$ and the set of morphisms from
$x$ to $y$ is given by $_yD_x=\sum_{g\in G}
\beta_g(j(_{g^{-1}y}C_{g^{-1}x}))$. It is clear that $D$ is a
small k-subcategory of $B$, $j:C\to D$ is a faithful functor and
$\beta$ is a global action of $G$ on $D$. Also condition (iv) of
the definition of enveloping action is fulfilled.

Recall that $\alpha^{h^{-1}}(_{y}\mathcal{I}_{x}^{h} \cap
{_{y}\mathcal{I}_{x}^{g^{-1}}})={_{h^{-1}y}\mathcal{I}_{h^{-1}x}^{(gh)^{-1}}}\cap
{_{h^{-1}y}\mathcal{I}_{h^{-1}x}^{h^{-1}}}$. Using this and taking
$h^{-1}=g$ we easily obtain $\alpha^g(1^{g^{-1}}1^{h})=1^g1^{gh}$,
for any $g,h\in G$. To simplify notation we will write
$\mathcal{I}^g$ instead of $_x\mathcal{I}^g_y$, for $x,y\in C_0$.

For any $\eta\in \mathcal{I}^{g^{-1}}$ we have $\alpha^g(\eta
1^{g^{-1}h})\in \alpha^g(\mathcal{I}^{g^{-1}}\cap
\mathcal{I}^{g^{-1}h})=\mathcal{I}^g\cap \mathcal{I}^h$. Thus
using conditions (v) and (vi) of Definition 3.1 we obtain
$$\alpha^{h^{-1}g}(\eta
1^{g^{-1}h})=\alpha^{h^{-1}}(\alpha^g(\eta
1^{g^{-1}h}))=\alpha^{h^{-1}}(\alpha^g(\eta)\alpha^g(1^{g^{-1}}1^{g^{-1}h}))=$$
$$\alpha^{h^{-1}}(\alpha^g({\eta})1^g1^h)=\alpha^{h^{-1}}(\alpha^g(\eta)1^h).$$

Now we show condition (iv) of Definition 4.2. In fact, let $\eta\in
\mathcal{I}^{g^{-1}}$. Then
$\beta_g(j(\eta))\mid_h=j(\eta)\mid_{g^{-1}h}=\alpha^{h^{-1}g}(\eta
1^{g^{-1}h})=\alpha^{h^{-1}}(\alpha^g(\eta)1^h)=j(\alpha^g(\eta))\mid_h$.
Consequently $\beta_g(j(\eta))=j(\alpha^g(\eta))$, for any
$\eta\in \mathcal{I}^{g^{-1}}$ and (iv) holds.

Let see now that condition (iii) of Definition 4.2 holds. Let $\sigma \in j(_yC_x)\cap
\beta_g(j(_{g^{-1}y}C_{g^{-1}x}))$. Then there are $\eta\in
{_yC_x}$ and $\mu\in {_{g^{-1}y}C_{g^{-1}x}}$ such that
$\sigma=j(\eta)=\beta_g(j(\mu))$. Thus for any $h\in G$ we have
$j(\eta)\mid_h=\beta_g(j(\mu))\mid_h$. This implies that
$\alpha^{h^{-1}}(\eta
1^h)=j(\mu)\mid_{g^{-1}h}=\alpha^{h^{-1}g}(\mu 1^{g^{-1}h})$.
Hence taking $h=e$ we obtain $\eta=\alpha^g(\mu 1^{g^{-1}})\in
\mathcal{I}^g$ and so $\sigma\in j(\mathcal{I}^g)$. The argument
shows that $j(C)\cap \beta_g(j(C))\subseteq j(\mathcal{I}^g)$.

On the other hand, if $\nu=j(\eta)$ for $\eta\in \mathcal{I}^g$ we have that
$\nu\in j(C)$. Also there exists $\tau\in \mathcal{I}^{g^{-1}}$
with $\alpha_g(\tau)=\eta$. Hence, $\nu=j(\alpha_g(\tau))$ and so
for any $h\in G$ we have
$\nu\mid_h=j(\alpha_g(\tau))\mid_h=\beta_g(j(\tau)\mid_h$.
So, $\nu\in \beta_g(j(C))$ and the relation (iii) follows.

Finally, we see that also condition (ii) is satisfied. Let $\eta\in
{_{g^{-1}z}C_{g^{-1}y}}$ and $\mu\in {_yC_x}$. Then
$$\beta_g(j(\eta))\mid_h\circ
j(\mu)\mid_h=j(\eta)\mid_{g^{-1}h}\circ j(\mu)\mid_h
=\alpha^{h^{-1}g}(\eta 1^{g^{-1}h})\alpha^{h^{-1}}(\mu 1^h)=$$ $$
\alpha^{h^{-1}}(\alpha^g(\eta 1^{g^{-1}}) 1^h)\alpha^{h^{-1}}(\mu
1^h)= \alpha^{h^{-1}}(\alpha^g(\eta 1^{g^{-1}}) \mu 1^h)=
j(\alpha^g(\eta 1^{g^{-1}})\mu)\mid_h.$$ \hspace{.5cm}
Consequently $\beta_g(j(\eta))\circ j(\mu)=j(\alpha^g(\eta
1^{g^{-1}})\mu)\in j(C)$. Similarly we prove $j(\mu) \circ
\beta_g(j(\eta)) =j(\mu\alpha^g(\eta 1^{g^{-1}}))\in j(C)$. Thus
condition (ii) also holds.

It remains to prove the uniqueness. Assume that $(D,\beta)$ and
$(E,\gamma)$ are two globalizations of $\alpha$, where $j:C\to D$
and $j^{'}:C\to E$ are the canonical functors of Definition 4.2.
Then $_{y}D_{x}=\sum_{g\in G} \beta_g(j(_{g^{-1}y}C_{g^{-1}x}))$ and
$_{y}E_{x}=\sum_{g\in G} \gamma_g(j^{'}(_{g^{-1}y}C_{g^{-1}x}))$,
for any $x,y\in C_0=D_0=E_0$. We define the mapping
$\Theta:{_y}D_x\to {_{y}E_x}$ by $\Theta(\sum_{1\leq i\leq
n}\beta_{g_i}(j(f_i)))=\sum_{1\leq i\leq
n}\gamma_{g_i}(j^{'}(f_i))$, for $f_i\in {_{g^{-1}y}}C_{g^{-1}x}$,
$1\leq i\leq n$. As in the last part of the proof of Theorem 4.5 of
\cite{DE} it follows that $\Theta$ is well-defined and so it is an
isomorphism of $k$-modules. This completes the proof. \sq

\vspace{.4cm}

We were unable either to prove or to disprove the result
corresponding to Theorem 4.5, in general. The implication in one
direction always holds. In fact, to end the paper we can easily to
obtain the following, repeating the first part of the proof of
Theorem 4.5.

\begin{prop} Let $\alpha$ be a partial action of a group $G$ on
$C$ and suppose that the partial action $\alpha$ has an enveloping
action $(D,\beta)$. Then for each $x\in C_0,
{_{x}\mathcal{I}^{g}_{x}}$ has a local identity, for all $g\in
G$.\end{prop}


\begin{thebibliography}{99}

\bibitem{A} Abadie, F., Enveloping actions and Takai duality for
partial actions, J. Funct. Analysis 197 (2003), 14-67.

\bibitem{ME} Alves, Marcelo M. S., Batista, E., Enveloping Actions for Partial Hopf Actions,
 Comm. in Algebra, 38 (8)(2010), 1532-4125, 2872–2902.
 
\bibitem{CJ} Caenepeel, S. and Janssen, K., Partial entwining
structures, Comm. in Algebra, 36(8) (2010), 2923–2946



\bibitem{CM} Cibilis, C. and Marcos, E.N., Skew category, Galois
covering and smash product of a k-category. Proccedings of
American Mathematical Society 1134(1)(2005), 39-50.

\bibitem{DEP} Dokuchaev, M., Exel, R. and Piccione, P., Partial
representations and partial group algebras, Journal of Algebra
226(2000), 251-258.

\bibitem{DE} Dokuchaev, M. and Exel, R., Associativity of crossed
products by partial actions, enveloping actions and partial
representations, Trans. Amer. Math. Soc. 357(2005), 1931-1952.

\bibitem{DFP} Dokuchaev, M., Ferrero, M. and Paques A., Partial
actions and Galois theory, Journal of Pure and applied Algebra, to
appear.

\bibitem{E} Exe, R., Twisted partial actions: a classification of
regular C*-algebraic bundles. Proceedings of London Mathematical
society 74 (1997), 417-443.

\bibitem{G} Gordon, R., G-categories, Memoirs of American
Mathematical society 482 (1993).

\bibitem{GM} Green, E. L., Marcos, E. N., Graded quotients of path algebras: A Local theory.
Journal of Pure an Applied Algebra 93 (1994) 195-226.


\end{thebibliography}
\end{document}